\documentclass{amsart}
\usepackage{amsmath}
\usepackage[all]{xy}

\theoremstyle{plain}
\newtheorem{theorem}{Theorem}
\newtheorem{lemma}{Lemma}
\newtheorem*{theoremTTS}{Theorem TTS}
\newtheorem*{theoremTTS'}{Theorem TTS$'$}
\newtheorem*{theoremSDS}{Theorem SDS}
\newtheorem*{theoremSDS'}{Theorem SDS$'$}
\newtheorem*{theoremGP}{Theorem GP}
\newtheorem*{theoremGP'}{Theorem GP$'$}
\newtheorem*{conjectureTTS}{Conjecture TTS}
\newtheorem*{conjectureTTS'}{Conjecture TTS$'$}

\theoremstyle{definition}
\newtheorem{definition}[lemma]{Definition}
\newtheorem{problem}{Problem}
\newtheorem*{remark}{Remark}

\begin{document}

\title{Twin Towers of Hanoi}

\author{Zoran {\v S}uni\'c}

\address{Department of Mathematics, Texas A\&M University, College Station, TX 77843-3368, USA}
\email{sunic@math.tamu.edu}

\thanks{This material is based upon work supported by the National Science Foundation}

\dedicatory{Dedicated to Antonio Machi on the occasion of his retirement}

\begin{abstract}
In the Twin Towers of Hanoi version of the well known Towers of Hanoi Problem
there are two coupled sets of pegs. In each move, one chooses a pair of pegs in
one of the sets and performs the only possible legal transfer of a disk
between the chosen pegs (the smallest disk from one of the pegs is moved to the
other peg), but also, simultaneously, between the corresponding pair of pegs in the coupled set (thus the same sequence of moves is always used in both sets). We provide upper and lower bounds on the length of the optimal solutions to problems of the following type. Given an initial and a final position 
of $n$ disks in each of the coupled sets, what is the smallest number of
moves needed to simultaneously obtain the final position from the initial one in each set? 
Our analysis is based on the use of a group, called Hanoi Towers group, of rooted ternary tree
automorphisms, which models the original problem in such a way that the
configurations on $n$ disks are the vertices at level $n$ of the tree and the
action of the generators of the group represents the three possible moves
between the three pegs. The twin version of the problem is analyzed by considering the action of Hanoi Towers group on pairs of vertices. 
\end{abstract}

\maketitle


\section{Towers of Hanoi and Twin Towers of Hanoi}

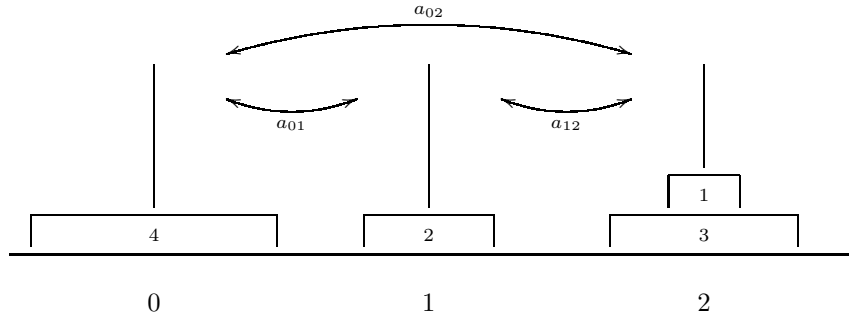
\begin{figure}[!ht]
\[
\xymatrix@C=5pt@R=9pt{
&&&&&\ar@{-}[dddd] &&\ar@{<->}@/^1pc/[rrrrrrrrrrrrrr]^{a_{02}}&&& &&&&\ar@{-}[dddd]& &&&&& &&&\ar@{-}[ddd]&& &&& \\
&&&&& &&\ar@{<->}@/_/[rrrrr]_{a_{01}}&&& &&&&& &\ar@{<->}@/_/[rrrrr]_{a_{12}}&&&& &&&&& &&& \\
&&&&& &&&&& &&&&& &&&&& &&&&& &&& \\
&&&&& &&&&& &&&&& &&&&& &&*{}\ar@{-}[rr]\ar@{-}[d]&&*{}\ar@{-}[d]& &&& \\
&*{} \ar@{-}[rrrrrrrr]\ar@{-}[d]&&&& &&&&*{}\ar@{-}[d]& &&*{}\ar@{-}[rrrr]\ar@{-}[d]&&& &*{}\ar@{-}[d]&&&& *{}\ar@{-}[rrrrrr]\ar@{-}[d]&&&\ar@{}[u]|{1}&& &*{}\ar@{-}[d]&& \\
\ar@{-}[rrrrrrrrrrrrrrrrrrrrrrrrrrrr]&&&&&\ar@{}[u]|{4} &&&&& &&&&\ar@{}[u]|{2}& &&&&& &&&\ar@{}[u]|{3}&& &&& \\
&&&&&0 &&&&& &&&&1& &&&&& &&&2&& &&&
}
\]
\caption{A configuration on four disks}
\label{f:configuration}
\end{figure}

We first describe the well known Hanoi Towers Problem on $n$ disks and 3 pegs.
The $n$ disks have different size. Allowed positions (which we call
configurations) of the disks on the pegs are those in which no disk is
on top of a smaller disk. An example of a configuration on 4 disks is provided
in Figure~\ref{f:configuration}). In a single move, the top disk from one of the
pegs can be transferred to the top position on another peg as long as the
newly obtained position of the disks is allowed (it is a configuration). 


Label the three pegs by 0, 1 and 2. At any moment, regardless of the current
configuration, there are exactly three possible moves, denoted by $a_{01}$,
$a_{02}$, and $a_{12}$. The move $a_{ij}$ transfers the smallest disk from pegs
$i$ and $j$ between these two pegs. More precisely, if the smallest disk on
pegs $i$ and $j$ is on $i$ the move $a_{ij}$ transfers it to $j$, and if it
is on $j$ the move transfers it to $i$. For instance, the move $a_{01}$
applied to the configuration in Figure~\ref{f:configuration} transfers disk 2
from peg 1 to peg 0, $a_{02}$ transfers disk 1 from peg 2 to peg 0, and
$a_{12}$
transfers disk 1 from peg 2 to peg 1. We do not need to specify the direction of
the transfer, since it is uniquely determined by the disks (by their size) that
are currently on pegs $i$ and $j$. In the exceptional case when there are no
disks on either peg $i$ or $j$, the move $a_{ij}$ leaves such a configuration
unchanged. 

In the classical Towers of Hanoi Problem on $n$ disks all disks are initially on
one of the pegs and the goal is to transfer all of them to another (prescribed) peg in the
smallest possible number of moves. It is well known that the optimal solution is
unique and consists of $2^n-1$ moves. One may pose a more general problem such
as, given some initial and final configurations on $n$ disks, what is the
smallest number of moves needed to obtain the final configuration from the
initial one. It turns out that this problem always has a solution (regardless of
the chosen initial and final configurations) and that the optimal solution is
either unique or there are exactly two solutions. The latter happens for a
relatively small number of choices of initial and final configurations. For a
survey on topics and results related to Hanoi Towers Problem see~\cite{hinz:ens}
and for an optimal solution (represented/obtained by a finite automaton) for any
pair of configurations see~\cite{romik:hanoi}. Note that, in this setting, none
of the instances of the general problem is more difficult (in terms of
the optimal number of moves) than the classical problem. 

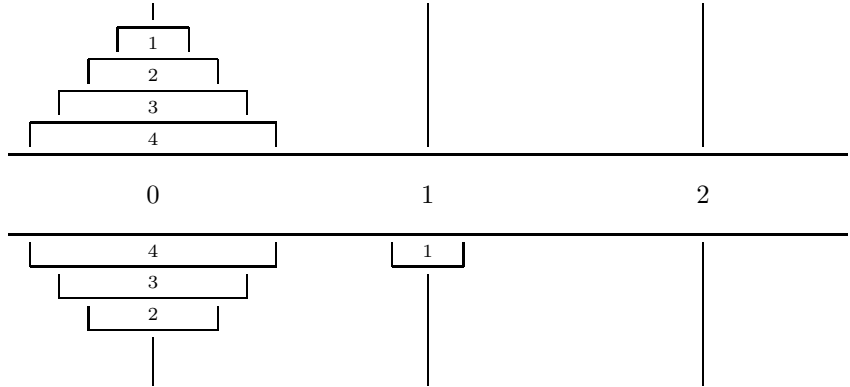
\begin{figure}[!ht]
\[
\xymatrix@C=5pt@R=6pt{
&&&&&\ar@{-}[d] &&&&& &&&&\ar@{-}[ddddd]& &&&&& &&&\ar@{-}[ddddd]&& &&& 
\\
&&&&*{}\ar@{-}[rr]\ar@{-}[d]&\ar@{}[d]|{1} &*{}\ar@{-}[d]&&&& &&&&& &&&&& &&&&& &&& 
\\
&&&*{}\ar@{-}[rrrr]\ar@{-}[d]&&\ar@{}[d]|{2} &&*{}\ar@{-}[d]&&& &&&&& &&&&& &&&&& &&& \\
&&*{}\ar@{-}[rrrrrr]\ar@{-}[d]&&&\ar@{}[d]|{3} &&&*{}\ar@{-}[d]&& &&&&& &&&&& &&&&& &&& \\
&*{}\ar@{-}[rrrrrrrr]\ar@{-}[d]&&&&\ar@{}[d]|{4} &&&&*{}\ar@{-}[d]& &&&&& &&&&& &&&&& &&& \\
\ar@{-}[rrrrrrrrrrrrrrrrrrrrrrrrrrrr]&&&&& &&&&& &&&&& &&&&& &&&&& &&& 
\\
&&&&&0 &&&&& &&&&1& &&&&& &&&2&& &&& \\
\ar@{-}[rrrrrrrrrrrrrrrrrrrrrrrrrrrr]&&&&&\ar@{}[d]|{4}&&&&& &&&&\ar@{}[d]|{1}& &&&&& &&&\ar@{-}[ddddd]&& &&& 
\\
&*{}\ar@{-}[rrrrrrrr]\ar@{-}[u]&&&&\ar@{}[d]|{3} &&&&*{}\ar@{-}[u]& &&&*{}\ar@{-}[rr]\ar@{-}[u]&\ar@{-}[dddd]&*{}\ar@{-}[u] &&&&& &&&&& &&& \\
&&*{}\ar@{-}[rrrrrr]\ar@{-}[u]&&&\ar@{}[d]|{2} &&&*{}\ar@{-}[u]&& &&&&& &&&&& &&&&& &&& \\
&&&*{}\ar@{-}[rrrr]\ar@{-}[u]&&\ar@{-}[dd] &&*{}\ar@{-}[u]&&& &&&&& &&&&& &&&&& &&& \\
&&&&& &&&&& &&&&& &&&&& &&&&& &&& \\
&&&&& &&&&& &&&&& &&&&& &&&&& &&& \\
}
\]
\caption{Initial position for the Small Disk Shift Problem}
\label{f:initial}
\end{figure}

In the Twin Towers of Hanoi version two sets of three pegs labeled by 0, 1 and 2
are coupled up. We often refer to the two sets as the top and the bottom set. A
coupled configuration on $n$ disks is a pair of configurations on $n$ disks, one
in each set (see, for instance, the coupled  configuration on 4 disks in
Figure~\ref{f:initial}). A move $a_{ij}$ applied to a coupled configuration
consists of application of the move $a_{ij}$ to each configuration in the
coupled pair. For instance, the move $a_{01}$ applied to the coupled
configuration in Figure~\ref{f:initial} transfers disk 1 in the top set to peg 1
and, simultaneously, disk 1 in the bottom set to peg 0. The move
$a_{02}$ applied to the same coupled configuration, transfers disk 1 in the top
set and disk 2 in the bottom set to peg 2 (in their sets), and $a_{12}$
changes nothing in the top set and transfers disk 1 in the bottom set to peg
2. 


In the setting of Twin Towers we pose three problems. 

\begin{problem}[Twin Towers Switch]
Given the initial coupled configuration in which all disks in the top set are on
peg 0 and all disks in the bottom set are on peg 2, how many moves are needed to
obtain the final coupled configuration in which all disks in the top set are on
peg 2, and all disks in the bottom set are on peg 0?  
\end{problem}

Note that the Twin Towers Switch Problem asks for simultaneous solution of two
instances of the classical Hanoi Towers Problem (all disks are, simultaneously,
using the same sequence of moves, transferred from peg 0 to peg 2 in the top
set, and from peg 2 to peg 0 in the bottom set). 

\begin{problem}[Small Disk Shift]
Given the initial coupled configuration in Figure~\ref{f:initial}, how many
moves are needed to obtain the final coupled configuration in which all disks
are in the same positions as in the initial one, except the smallest disk in
each set is moved one peg to the right (disk 1 in the top configuration to peg
1, and disk 1 in the bottom configuration to peg 2)? 
\end{problem}

\begin{problem}[General Problem]
Given any initial coupled configuration and any final coupled configuration what
is the smallest number of moves needed to obtain the final configuration from
the initial one?  
\end{problem}

We provide an upper bound for the Twin Towers Switch, exact answer for the Small
Disk Shift, and lower and upper bounds for the General Problem restricted to basic coupled configurations (defined below).  

\begin{theoremTTS}[Twin Towers Switch]
The smallest number of moves needed to solve the Twin Towers Switch Problem on
$n$ disks is no greater than $a(n)$, where  
\[
a(n) =  
 \begin{cases}
   1, & n=1, \\
   \frac{4}{3} \cdot 2^n -\frac{(-1)^n}{3} , & n \geq 2.
 \end{cases} 
\]
\end{theoremTTS}

\begin{remark}
The sequence $a(n)$ satisfies the Jacobshtal linear recursion
\[
 a(n) = a(n-1) + 2a(n-2), \qquad \text{ for } n \geq 4,  
\]
with initial condition $a(1)=1$, $a(2)=5$, and $a(3)=11$.
\end{remark}

\begin{conjectureTTS}
The smallest number of moves needed to solve the Twin Towers Switch problem on
$n$ disks is exactly $a(n)$.  
\end{conjectureTTS}

Note that the Twin Towers Switch, requiring no more than roughly $\frac{4}{3}2^n$ moves is
not considerably more difficult than the classical problem of moving a single
tower, which requires roughly $2^n$ moves. In fact, there are more difficult
problems that can be posed in the context of coupled sets (recall that there
are no problems that are more difficult than the classical problem when only one
set of disks is considered). For instance, the next result implies that the
Small Disk Shift Problem requires more moves than the Twin Towers Switch Problem.  

\begin{theoremSDS}[Small Disk Shift]
The smallest number of moves $d(n)$ needed to solve the Small Disk Shift Problem
on $n$ disks is equal to 
\[
d(n) = 
 \begin{cases}
   2, & n=1, \\
   6, & n=2, \\
   2 \cdot 2^n, & n \geq 3.
 \end{cases} 
\]
\end{theoremSDS}

In order to state our result on the General Problem, we need the notion of
compatible coupled configurations. An initial coupled configuration $I$ on $n$
disks is \emph{compatible} to the final coupled configuration $F$ on $n$ disks if  $F$
can be obtained from $I$ in a finite number of moves. A coupled configuration
is called \emph{basic} if the smallest disks in its top configuration and the
smallest disk in its bottom configuration are not on corresponding pegs (it is
not the case that both are on peg 0, both on peg 1, or both on peg 2). 

Note that, based on the branching structure of Hanoi Towers group described by
Grigorchuk and the author in~\cite{grigorchuk-s:standrews}, D'Angeli and Donno
show in~\cite{d'angeli-d:hanoi} that Hanoi Towers group acts distance 2-transitively on the levels of the rooted ternary tree. This provides a 
characterization of the pairs of compatible coupled configurations. In
particular, their result implies that all basic coupled configurations are
compatible. We quote their result in more detail (Theorem~\ref{t:d'angeli-d}),
after we sufficiently develop the necessary terminology.  Along the way we
provide a different proof (we need it for our upper bound estimate on the
General Problem). Note that an interesting consequence of the result of D'Angeli
and Donno is that the Hanoi Towers group induces an infinite sequence of finite
Gel$'$fand pairs (see~\cite{d'angeli-d:hanoi} for details). 

\begin{theoremGP}[General problem for basic configurations]
The number of moves needed to obtain one basic coupled configuration on $n$
disks from another is no greater than  
\[
 \frac{11}{3} \times 2^n =  3.\overline{66} \times 2^n. 
\]
\end{theoremGP}

Note that the coupled configurations in Theorem~SDS are basic. Thus, Theorem~SDS implies that for at least one pair of basic coupled configurations the smallest number of moves that is needed to obtain one from the other is exactly $2 \times 2^n$. 

Obtaining good upper bound seems to be a difficult task, since one needs to solve all instances of the problem in optimal or nearly optimal way. Lower bounds seem a bit easier to obtain since they may be derived from lower bounds from some specific, well chosen, instances. The lower bound ($2 \times 2^n$) and the upper bound ($3.\overline{66} \times 2^n$) provided here differ by less than a factor of two. 

All results mentioned so far will be recast in the following sections in the natural setting of group actions on rooted trees. The reason is that this setting provides a convenient language and tools to prove our results. 

\subsection*{Acknowledgment}

The author is thankful to Tullio Ceccherini-Silberstein and Alfredo Donno for their help, useful remarks, and corrections. 


\section{Encoding by words and tree automorphisms}

We start by an encoding of the original Hanoi Towers Problem on three pegs, as
originally presented in~\cite{grigorchuk-s:hanoi-cr} (and further elaborated
in~\cite{grigorchuk-s:standrews,grigorchuk-s:hanoi-spectrum}), by a group of
rooted ternary tree automorphisms.  

Label the disks by $1,2,\dots,n$ according to their size (smallest to
largest). The configurations can be encoded by words over the finite alphabet
$X=\{0,1,2\}$. The letters in this alphabet represent the pegs. The word
$x_1x_2 \dots x_n$ represents the unique configuration on $n$ disks in which,
for $i=1,\dots,n$, the disk $i$ is on peg $x_i$. For example, the word 2120
represents the configuration in Figure~\ref{f:configuration}. Note that there
are exactly $3^n$ configurations on $n$ disks.  

The moves $a_{ij}$ are encoded as the transformations of the set of all finite
words $X^*$ over $X$ defined by  
\begin{alignat*}{6}
 &a_{01}(2 \dots 2 0u) &&=  2 \dots 2 1u, & \quad
 &a_{02}(1 \dots 1 0u) &&=  1 \dots 1 2u, & \quad
 &a_{12}(0 \dots 0 1u) &&=  0 \dots 0 2u,
 \\
 &a_{01}(2 \dots 2 1u) &&=  2 \dots 2 0u, & 
 &a_{02}(1 \dots 1 2u) &&=  1 \dots 1 0u, &
 &a_{12}(0 \dots 0 2u) &&=  0 \dots 0 1u, 
 \\
 &a_{01}(2 \dots 2) &&=  2 \dots 2, &
 &a_{02}(1 \dots 1) &&=  1 \dots 1, &
 &a_{12}(0 \dots 0) &&=  0 \dots 0, 
\end{alignat*}
for any word $u$ in $X^*$. Thus, $a_{ij}$ changes the first occurrence of $i$ or
$j$ to the other of these two symbols. The point of, say, $a_{01}$ ``ignoring''
initial prefixes of the form $2^\ell$ is that such prefixes represent small
disks on peg $2$, and $a_{01}$ should ignore such disks, since it is supposed to
transfer a disk between peg $0$ and peg $1$. The first occurrence of $0$ or $1$
represents the smallest disk on one of these two pegs and changing this
occurrence of the symbol 0 or 1 to the other one in the code of the given
configuration transfers the corresponding disk to the other peg. Note that if
$a_{ij}$ is applied to (a code of) a configuration that has no occurrences of
$i$ or $j$ it leaves such a configuration unchanged. This corresponds to
the situation in which there are no disks on pegs $i$ ad $j$ and the move $a_{ij}$
has no effect on such a configuration since there are no disks to be moved. 

In order to work with more compact notation, set  
\[
 a_{01} = a, \qquad a_{02} = b, \qquad a_{12} = c. 
\]
In this notation, the moves $a$, $b$ and $c$ act on the set of all finite words $X^*$ by 
\begin{alignat}{6}\label{e:formula}
 &a(2 \dots 2 0u) &&=  2 \dots 2 1u, & \quad
 &b(1 \dots 1 0u) &&=  1 \dots 1 2u, & \quad
 &c(0 \dots 0 1u) &&=  0 \dots 0 2u, \notag
 \\
 &a(2 \dots 2 1u) &&=  2 \dots 2 0u, & 
 &b(1 \dots 1 2u) &&=  1 \dots 1 0u, &
 &c(0 \dots 0 2u) &&=  0 \dots 0 1u, 
 \\
 &a(2 \dots 2) &&=  2 \dots 2, &
 &b(1 \dots 1) &&=  1 \dots 1, &
 &c(0 \dots 0) &&=  0 \dots 0. \notag
\end{alignat}

\begin{figure}[!ht]
\[
\xymatrix@C=15pt{
&&&&&&&
*[o][F-]{\bullet} \ar@{-}[dl]_{a} \ar@{-}[dr]^{c} \ar@{-}@(ul,ur)^{b} \ar@{}[r]|>{111}&&&&&&&
\\
&&&&&&
*[o][F-]{\bullet} \ar@{-}[dl]_{c} \ar@{-}[rr]_{b} \ar@{}[l]|>{011} && 
*[o][F-]{\bullet} \ar@{-}[dr]^{a} \ar@{}[r]|>{211}&&&&&&
\\
&&&&&
 *[o][F-]{\bullet} \ar@{-}[dl]_{b} \ar@{-}[dr]^{a} \ar@{}[l]|>{021} &&&&
 *[o][F-]{\bullet} \ar@{-}[dl]_{c} \ar@{-}[dr]^{b} \ar@{}[r]|>{201} &&&&&
\\
&&&&
*[o][F-]{\bullet} \ar@{-}[dl]_{a} \ar@{-}[rr]_{c} \ar@{}[l]|>{221} && 
*[o][F-]{\bullet} \ar@{-}[rr]_{b} \ar@{}[d]|{121} &&
*[o][F-]{\bullet} \ar@{-}[rr]_{a} \ar@{}[d]|{101} && 
*[o][F-]{\bullet} \ar@{-}[dr]^{c} \ar@{}[r]|>{001}&&&&
\\
&&&
*[o][F-]{\bullet} \ar@{-}[dl]_{c} \ar@{-}[dr]^{b} \ar@{}[l]|>{220} &&&&&&&&
*[o][F-]{\bullet} \ar@{-}[dl]_{b} \ar@{-}[dr]^{a} \ar@{}[r]|>{002} &&&
\\
&&
*[o][F-]{\bullet} \ar@{-}[dl]_{b} \ar@{-}[rr]_{a} \ar@{}[l]|>{120} && 
*[o][F-]{\bullet} \ar@{-}[dr]^{c} \ar@{}[r]|>{020} &&&&&&
*[o][F-]{\bullet} \ar@{-}[dl]_{a} \ar@{-}[rr]_{c} \ar@{}[l]|>{202} && 
*[o][F-]{\bullet} \ar@{-}[dr]^{b} \ar@{}[r]|>{102}&&
 \\
 & 
 *[o][F-]{\bullet} \ar@{-}[dl]_{a} \ar@{-}[dr]^{c} \ar@{}[l]|>{100} &&&&
 *[o][F-]{\bullet} \ar@{-}[dl]_{b} \ar@{-}[dr]^{a} \ar@{}[r]|>{010} &&&& 
 *[o][F-]{\bullet} \ar@{-}[dl]_{c} \ar@{-}[dr]^{b} \ar@{}[l]|>{212} &&&&
 *[o][F-]{\bullet} \ar@{-}[dl]_{a} \ar@{-}[dr]^{c} \ar@{}[r]|>{122} &
 \\ 
 *[o][F-]{\bullet} \ar@{-}[rr]_{b} \ar@{-}@(ul,dl)_{c} \ar@{}[d]|{000} && 
 *[o][F-]{\bullet} \ar@{-}[rr]_{a} \ar@{}[d]|{200} &&
 *[o][F-]{\bullet} \ar@{-}[rr]_{c} \ar@{}[d]|{210} && 
 *[o][F-]{\bullet} \ar@{-}[rr]_{b} \ar@{}[d]|{110} &&  
 *[o][F-]{\bullet} \ar@{-}[rr]_{a} \ar@{}[d]|{112} && 
 *[o][F-]{\bullet} \ar@{-}[rr]_{c} \ar@{}[d]|{012} &&
 *[o][F-]{\bullet} \ar@{-}[rr]_{b} \ar@{}[d]|{022} && 
 *[o][F-]{\bullet} \ar@{-}@(ur,dr)^{a} \ar@{}[d]|{222}
 \\
 &&&&& &&&&& &&&&
}
\]
\caption{$\Gamma_3$, the Hanoi graph on 3 disks}
\label{f:schreier3}
\end{figure}
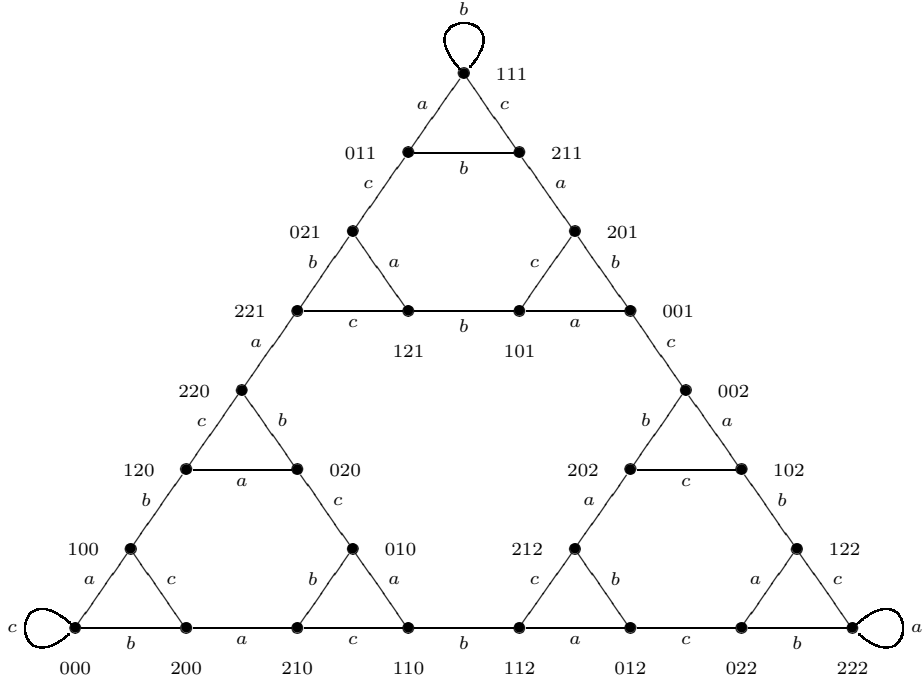

Hanoi graph on $n$ disks, denoted by $\Gamma_n$, is the graph on
$3^n$ vertices representing the configurations on $n$ disks. Two vertices $u$
and $v$ are connected by an edge labeled by $s \in \{a,b,c\}$ if the
configurations represented by $u$ and $v$ can be obtained from each other by
application of the move $s$ (note that each of the moves is an involution). The
Hanoi graph on 3 disks is depicted in Figure~\ref{f:schreier3}. 
Graphs very similar to the graphs we just defined have already appeared in the
literature in connection to Hanoi Towers Problem (see, for instance,
\cite{hinz:ens}). The difference is that the edges are usually not labeled and there are no loops at the corners. 

The set of all words $X^*$ has the structure of a rooted ternary tree in which
the root is the empty word, level $n$ of the tree consists of the $3^n$ words of
length $n$ over $X$, and each vertex (each word) $u$ has three children, $u0$,
$u1$ and $u2$. The transformations $a$, $b$ and $c$ act on the tree $X^*$ as
tree automorphisms (in particular, they preserve the root and the levels of the
tree). Thus, $a$, $b$ and $c$ generate a group of automorphisms of the rooted
ternary tree $X^*$. The group $H=\langle a,b,c \rangle$, called Hanoi Towers group,
was defined in~\cite{grigorchuk-s:hanoi-cr}. The Hanoi graph $\Gamma_n$ is the
Schreier graph, with respect to the generating set $\{a,b,c\}$, of the action of
$H$ on the words of length $n$ in $X^*$ (Schreier graph of the action on level
$n$ in the tree).  

A sequence of moves is a word over $S=\{a,b,c\}$. The order in which moves are
applied is from right to left as in the following calculation 
\[
 caba(0220) = cab(1220) = ca(1020) = c(0020) = 0010. 
\]

\begin{figure}[!ht]
\[
\xymatrix@C=3pt@R=6pt{
&&&&&&&&&&&&&&&&&
*[o][F-]{\bullet} \ar@{-}[dl]_{a} \ar@{-}[dr]^{c} \ar@{-}@(ul,ur)^{b}
\\
&&&&&&&&&&&&&&&&
*[o][F-]{\bullet} \ar@{-}[dl]_{c} \ar@{-}[rr]_{b} && *[o][F-]{\bullet} \ar@{-}[dr]^{a} 
\\
&&&&&&&&&&&&&&&
 *[o][F-]{\bullet} \ar@{-}[dl]_{b} \ar@{-}[dr]^{a} &&&&
 *[o][F-]{\bullet} \ar@{-}[dl]_{c} \ar@{-}[dr]^{b} &&&&&&&&&&&&& 
\\
&&&&&&&&&&&&&&
*[o][F-]{\bullet} \ar@{-}[dl]_{a} \ar@{-}[rr]_{c} && *[o][F-]{\bullet} \ar@{-}[rr]_{b} &&
*[o][F-]{\bullet} \ar@{-}[rr]_{a} && *[o][F-]{\bullet} \ar@{..}[dr]^{c} &&&&&&&&&&
\\
&&&&&&&&&&&&&&&&&&&&&&&&&&&&&&
\\
&&&&&&&&&&&&\ar@{-}[dl]_{c} \ar@{..}[ur]&&&&&&&&&&&&&&&&&&
\\
&&&&&&&&&&&
*[o][F-]{\bullet} \ar@{-}[dl]_{b} \ar@{-}[dr]^{a} &&&&&&&&
 &&&&&&&&&&&
\\ 
&&&&&&&&&&
*[o][F-]{\bullet} \ar@{-}[dl]_{a} \ar@{-}[rr]_{c} && *[o][F-]{\bullet} \ar@{-}[dr]^{b} &&&&&&
 && &&&&&&&&&&
\\
&&&&&&&&&
 *[o][F-]{\bullet} \ar@{-}[dl]_{c} \ar@{-}[dr]^{b} &&&&
 *[o][F-]{\bullet} \ar@{-}[dl]_{a} \ar@{-}[dr]^{c} &&&& &&&& &&&&&&&&&
\\
&&&&&&&&
*[o][F-]{\bullet} \ar@{-}[dl]_{b} \ar@{-}[rr]_{a} && *[o][F-]{\bullet} \ar@{-}[rr]_{c} &&
*[o][F-]{\bullet} \ar@{-}[rr]_{b} && *[o][F-]{\bullet} \ar@{..}[rr]_{a} &&
&&&&&&&&&&&&&&
\\
&&&&&&&
*[o][F-]{\bullet} \ar@{-}[dl]_{a} \ar@{-}[dr]^{c} &&&&&&&&&&&&&&&&&&&&&&&
\\
&&&&&&
*[o][F-]{\bullet} \ar@{-}[dl]_{c} \ar@{-}[rr]_{b} && *[o][F-]{\bullet} \ar@{-}[dr]^{a} 
&&&&&&&&&&&&&&&&&&&&&&
\\
&&&&&
 *[o][F-]{\bullet} \ar@{-}[dl]_{b} \ar@{-}[dr]^{a} &&&&
 *[o][F-]{\bullet} \ar@{-}[dl]_{c} \ar@{-}[dr]^{b}  &&&&&&&&&&&&&&&&
\\
&&&&
*[o][F-]{\bullet} \ar@{-}[dl]_{a} \ar@{-}[rr]_{c} && *[o][F-]{\bullet} \ar@{-}[rr]_{b} &&
*[o][F-]{\bullet} \ar@{-}[rr]_{a} && *[o][F-]{\bullet} \ar@{-}[dr]^{c} &&&&&&&&&&
\ar@{..}[dl]_{b} && &&
&&&&&& \ar@{..}[dr]^{c} &&&& 
\\
&&&
*[o][F-]{\bullet} \ar@{-}[dl]_{c} \ar@{-}[dr]^{b} &&&&&&&&
*[o][F-]{\bullet} \ar@{-}[dl]_{b} \ar@{-}[dr]^{a} &&&&&&&&
*[o][F-]{\bullet} \ar@{-}[dl]_{a} \ar@{-}[dr]^{c} &&&&&&&&&&&&
*[o][F-]{\bullet} \ar@{-}[dl]_{b} \ar@{-}[dr]^{a} &&&
\\
&&
*[o][F-]{\bullet} \ar@{-}[dl]_{b} \ar@{-}[rr]_{a} && *[o][F-]{\bullet} \ar@{-}[dr]^{c} &&&&&&
*[o][F-]{\bullet} \ar@{-}[dl]_{a} \ar@{-}[rr]_{c} && *[o][F-]{\bullet} \ar@{-}[dr]^{b} &&&&&&
*[o][F-]{\bullet} \ar@{-}[dl]_{c} \ar@{-}[rr]_{b} && *[o][F-]{\bullet} \ar@{-}[dr]^{a} &&&&&&&&&&
*[o][F-]{\bullet} \ar@{-}[dl]_{a} \ar@{-}[rr]_{c} && *[o][F-]{\bullet} \ar@{-}[dr]^{b} &&
 \\
 & 
 *[o][F-]{\bullet} \ar@{-}[dl]_{a} \ar@{-}[dr]^{c} &&&&
 *[o][F-]{\bullet} \ar@{-}[dl]_{b} \ar@{-}[dr]^{a} &&&& 
 *[o][F-]{\bullet} \ar@{-}[dl]_{c} \ar@{-}[dr]^{b} &&&&
 *[o][F-]{\bullet} \ar@{-}[dl]_{a} \ar@{-}[dr]^{c} &&&&
 *[o][F-]{\bullet} \ar@{-}[dl]_{b} \ar@{-}[dr]^{a} &&&& 
 *[o][F-]{\bullet} \ar@{-}[dl]_{c} \ar@{-}[dr]^{b} &&&&&&&&
  *[o][F-]{\bullet} \ar@{-}[dl]_{c} \ar@{-}[dr]^{b} &&&&
 *[o][F-]{\bullet} \ar@{-}[dl]_{a} \ar@{-}[dr]^{c} &
 \\ 
 *[o][F-]{\bullet} \ar@{-}[rr]_{b} \ar@{-}@(l,d)_{c} && *[o][F-]{\bullet} \ar@{-}[rr]_{a} &&
 *[o][F-]{\bullet} \ar@{-}[rr]_{c} && *[o][F-]{\bullet} \ar@{-}[rr]_{b} &&  
 *[o][F-]{\bullet} \ar@{-}[rr]_{a} && *[o][F-]{\bullet} \ar@{-}[rr]_{c} &&
 *[o][F-]{\bullet} \ar@{-}[rr]_{b}&& *[o][F-]{\bullet} \ar@{-}[rr]_{a} &&
 *[o][F-]{\bullet} \ar@{-}[rr]_{c} && *[o][F-]{\bullet} \ar@{-}[rr]_{b} &&  
 *[o][F-]{\bullet} \ar@{-}[rr]_{a} && *[o][F-]{\bullet} \ar@{-}[rr]_{c} &&
 &&\ar@{..}[ll] \ar@{-}[rr]_{b}&&
 *[o][F-]{\bullet} \ar@{-}[rr]_{a}&& *[o][F-]{\bullet} \ar@{-}[rr]_{c} &&
 *[o][F-]{\bullet} \ar@{-}[rr]_{b} && *[o][F-]{\bullet} \ar@{-}@(r,d)^{a}
}
\]
\caption{Hanoi graph on odd number of disks}
\label{f:odd-general}
\end{figure}
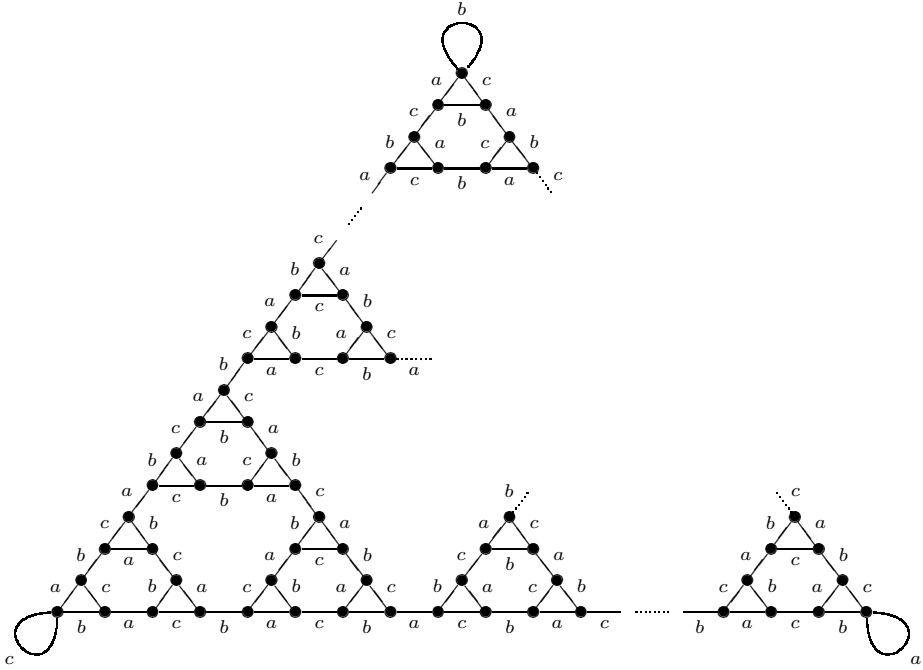

The structure of the Hanoi graphs is fairly well understood. In particular, for
$n \geq 0$, the Hanoi graph $\Gamma_{n+1}$ is obtained from the Hanoi graph
$\Gamma_{n}$ as follows~\cite{grigorchuk-s:standrews}. Three copies of $\Gamma_{n}$ are constructed by
appending the label 0, 1, and 2, respectively, to every vertex label in
$\Gamma_{n}$. Then the two loops labeled by $c$ at the vertices $0^{n}1$ and
$0^{n}2$ are deleted and replaced by an edge between $0^{n}1$ and $0^{n}2$
labeled by $c$, the two loops labeled by $b$ at the vertices $1^{n}0$ and
$1^{n}2$ are deleted and replaced by an edge between $1^{n}0$ and $1^{n}2$
labeled by $b$, and the two loops labeled by $a$ at the vertices $2^{n}0$ and
$2^{n}1$ are deleted and replaced by an edge between $2^{n}0$ and $2^{n}1$
labeled by $a$. Indeed, this ``rewiring'' on the next level (level $n+1$) needs to be done as
indicated since $c(0^{n}1)=0^{n}2$, $b(1^{n}0)= 1^{n}2$ and $a(2^{n}0)=2^{n}1$.
In general, the graphs for even and odd $n$ have the form provided in
Figure~\ref{f:even-general} and Figure~\ref{f:odd-general}.  These 
figures suffice for our purposes, since only the region near the path from $0^n$ to
$2^n$ (near the bottom) and near the path from $0^n$ to $1^n$ (near the left
side) play significant role in our considerations.  

\begin{figure}[!ht]
\[
\xymatrix@C=3pt@R=6pt{
&&&&&&&&&&&&&&&&&
*[o][F-]{\bullet} \ar@{-}[dl]_{c} \ar@{-}[dr]^{a} \ar@{-}@(ul,ur)^{b}
\\
&&&&&&&&&&&&&&&&
*[o][F-]{\bullet} \ar@{-}[dl]_{a} \ar@{-}[rr]_{b} && *[o][F-]{\bullet} \ar@{-}[dr]^{c} 
\\
&&&&&&&&&&&&&&&
 *[o][F-]{\bullet} \ar@{-}[dl]_{b} \ar@{-}[dr]^{c} &&&&
 *[o][F-]{\bullet} \ar@{-}[dl]_{a} \ar@{-}[dr]^{b} &&&&&&&&&&&&& 
\\
&&&&&&&&&&&&&&
*[o][F-]{\bullet} \ar@{-}[dl]_{c} \ar@{-}[rr]_{a} && *[o][F-]{\bullet} \ar@{-}[rr]_{b} &&
*[o][F-]{\bullet} \ar@{-}[rr]_{c} && *[o][F-]{\bullet} \ar@{..}[dr]^{a} &&&&&&&&&&
\\
&&&&&&&&&&&&&&&&&&&&&&&&&&&&&&
\\
&&&&&&&&&&&&\ar@{-}[dl]_{c} \ar@{..}[ur]&&&&&&&&&&&&&&&&&&
\\
&&&&&&&&&&&
*[o][F-]{\bullet} \ar@{-}[dl]_{a} \ar@{-}[dr]^{b} &&&&&&&&
 &&&&&&&&&&&
\\ 
&&&&&&&&&&
*[o][F-]{\bullet} \ar@{-}[dl]_{b} \ar@{-}[rr]_{c} && *[o][F-]{\bullet} \ar@{-}[dr]^{a} &&&&&&
 && &&&&&&&&&&
\\
&&&&&&&&&
 *[o][F-]{\bullet} \ar@{-}[dl]_{c} \ar@{-}[dr]^{a} &&&&
 *[o][F-]{\bullet} \ar@{-}[dl]_{b} \ar@{-}[dr]^{c} &&&& &&&& &&&&&&&&&
\\
&&&&&&&&
*[o][F-]{\bullet} \ar@{-}[dl]_{a} \ar@{-}[rr]_{b} && *[o][F-]{\bullet} \ar@{-}[rr]_{c} &&
*[o][F-]{\bullet} \ar@{-}[rr]_{a} && *[o][F-]{\bullet} \ar@{..}[rr]_{b} &&
&&&&&&&&&&&&&&
\\
&&&&&&&
*[o][F-]{\bullet} \ar@{-}[dl]_{b} \ar@{-}[dr]^{c} &&&&&&&&&&&&&&&&&&&&&&&
\\
&&&&&&
*[o][F-]{\bullet} \ar@{-}[dl]_{c} \ar@{-}[rr]_{a} && *[o][F-]{\bullet} \ar@{-}[dr]^{b} 
&&&&&&&&&&&&&&&&&&&&&&
\\
&&&&&
 *[o][F-]{\bullet} \ar@{-}[dl]_{a} \ar@{-}[dr]^{b} &&&&
 *[o][F-]{\bullet} \ar@{-}[dl]_{c} \ar@{-}[dr]^{a}  &&&&&&&&&&&&&&&&
\\
&&&&
*[o][F-]{\bullet} \ar@{-}[dl]_{b} \ar@{-}[rr]_{c} && *[o][F-]{\bullet} \ar@{-}[rr]_{a} &&
*[o][F-]{\bullet} \ar@{-}[rr]_{b} && *[o][F-]{\bullet} \ar@{-}[dr]^{c}  &&&&&&&&&&
\ar@{..}[dl]_{a} && &&
&&&&&& \ar@{..}[dr]^{b} &&&& 
\\
&&&
*[o][F-]{\bullet} \ar@{-}[dl]_{c} \ar@{-}[dr]^{a} &&&&&&&&
*[o][F-]{\bullet} \ar@{-}[dl]_{a} \ar@{-}[dr]^{b} &&&&&&&&
*[o][F-]{\bullet} \ar@{-}[dl]_{b} \ar@{-}[dr]^{c} &&&&&&&&&&&&
*[o][F-]{\bullet} \ar@{-}[dl]_{c} \ar@{-}[dr]^{a} &&&
\\
&&
*[o][F-]{\bullet} \ar@{-}[dl]_{a} \ar@{-}[rr]_{b} && *[o][F-]{\bullet} \ar@{-}[dr]^{c} &&&&&&
*[o][F-]{\bullet} \ar@{-}[dl]_{b} \ar@{-}[rr]_{c} && *[o][F-]{\bullet} \ar@{-}[dr]^{a} &&&&&&
*[o][F-]{\bullet} \ar@{-}[dl]_{c} \ar@{-}[rr]_{a} && *[o][F-]{\bullet} \ar@{-}[dr]^{b} &&&&&&&&&&
*[o][F-]{\bullet} \ar@{-}[dl]_{a} \ar@{-}[rr]_{b} && *[o][F-]{\bullet} \ar@{-}[dr]^{c} &&
 \\
 & 
 *[o][F-]{\bullet} \ar@{-}[dl]_{b} \ar@{-}[dr]^{c} &&&&
 *[o][F-]{\bullet} \ar@{-}[dl]_{a} \ar@{-}[dr]^{b} &&&& 
 *[o][F-]{\bullet} \ar@{-}[dl]_{c} \ar@{-}[dr]^{a} &&&&
 *[o][F-]{\bullet} \ar@{-}[dl]_{b} \ar@{-}[dr]^{c} &&&&
 *[o][F-]{\bullet} \ar@{-}[dl]_{a} \ar@{-}[dr]^{b} &&&& 
 *[o][F-]{\bullet} \ar@{-}[dl]_{c} \ar@{-}[dr]^{a} &&&&&&&&
  *[o][F-]{\bullet} \ar@{-}[dl]_{b} \ar@{-}[dr]^{c} &&&&
 *[o][F-]{\bullet} \ar@{-}[dl]_{a} \ar@{-}[dr]^{b} &
 \\ 
 *[o][F-]{\bullet} \ar@{-}[rr]_{a} \ar@{-}@(l,d)_{c} && *[o][F-]{\bullet} \ar@{-}[rr]_{b} &&
 *[o][F-]{\bullet} \ar@{-}[rr]_{c} && *[o][F-]{\bullet} \ar@{-}[rr]_{a} &&  
 *[o][F-]{\bullet} \ar@{-}[rr]_{b} && *[o][F-]{\bullet} \ar@{-}[rr]_{c} &&
 *[o][F-]{\bullet} \ar@{-}[rr]_{a}&& *[o][F-]{\bullet} \ar@{-}[rr]_{b} &&
 *[o][F-]{\bullet} \ar@{-}[rr]_{c} && *[o][F-]{\bullet} \ar@{-}[rr]_{a} &&  
 *[o][F-]{\bullet} \ar@{-}[rr]_{b} && *[o][F-]{\bullet} \ar@{-}[rr]_{c} &&
 &&\ar@{..}[ll] \ar@{-}[rr]_{c}&&
 *[o][F-]{\bullet} \ar@{-}[rr]_{a}&& *[o][F-]{\bullet} \ar@{-}[rr]_{b} &&
 *[o][F-]{\bullet} \ar@{-}[rr]_{c} && *[o][F-]{\bullet} \ar@{-}@(r,d)^{a}
}
\]
\caption{Hanoi graph on even number of disks}
\label{f:even-general}
\end{figure}
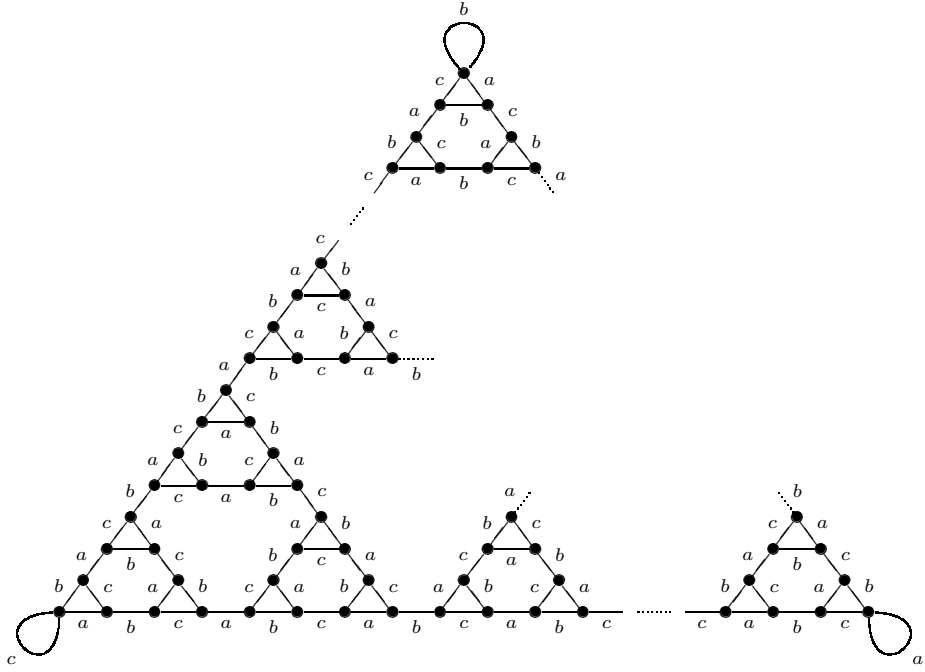

The following lemma, providing a non-recursive, optimal solution to the
classical Hanoi Towers Problem is part of the folklore (it has been proved and expressed in
many disguises and our setting may be considered one of them).  

\begin{lemma}\label{l:classic}
The diameter of the Hanoi Towers graph $\Gamma_n$ is $2^n-1$. It is achieved as the distance between any two of the configurations $0^n$, $1^n$, and
$2^n$. The unique sequence of moves of length $2^n-1$ between any two of these
configurations is given in the following table.  
\[
\begin{array}{c|ccc|ccc}
 & & \textup{even }n& & & \textup{odd }n & 
\\
 \textup{from} \backslash \textup{to} & 0^n & 1^n & 2^n & 0^n & 1^n & 2^n  
\\
  \hline 
 0^n & 
 \times & (cab)^{m(n)} & (cba)^{m(n)} &  
 \times & a(cba)^{m(n)} & b(cab)^{m(n)}
\\
 1^n & 
 (bac)^{m(n)} & \times & (bca)^{m(n)} &
 a(bca)^{m(n)} & \times & c(bac)^{m(n)} 
\\
 2^n & 
 (abc)^{m(n)} & (acb)^{m(n)} & \times & 
 b(acb)^{m(n)} & c(abc)^{m(n)} & \times 
 \\
 \hline
\end{array} 
\]
where $m(n)=\frac{1}{3}(2^n-1)$, for even $n$, and $m(n)=\frac{1}{3}(2^n-2)$, for odd $n$.  
\end{lemma}

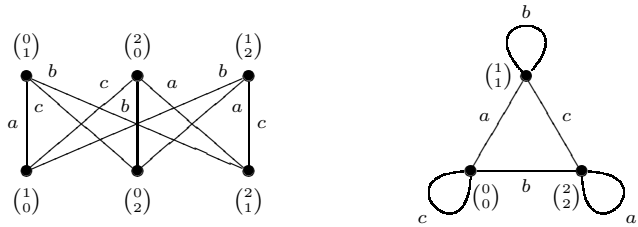
\begin{figure}
\[
\xymatrix@C=15pt@R=30pt{
 \ar@{}[d]|>>>>{\binom{0}{1}} && \ar@{}[d]|>>>>{\binom{2}{0}} && \ar@{}[d]|>>>>{\binom{1}{2}} &&&& && 
 \\
 *[o][F-]{\bullet} \ar@{-}[d]_{a} \ar@{-}[drr]_<<<<{c} \ar@{-}[drrrr]^<<<{b} && 
 *[o][F-]{\bullet} \ar@{-}[dll]_<<<<{c} \ar@{-}[d]_<<<<{b} \ar@{-}[drr]^<<<<{a} &&
 *[o][F-]{\bullet} \ar@{-}[dllll]_<<<{b} \ar@{-}[dll]^<<<<{a} \ar@{-}[d]^{c} &&&& & 
 *[o][F-]{\bullet} \ar@{-}[dl]_{a} \ar@{-}[dr]^{c} \ar@{-}@(ul,ur)^{b} \ar@{}[l]|{\binom{1}{1}} & 
 \\
 *[o][F-]{\bullet} && *[o][F-]{\bullet} && *[o][F-]{\bullet}  &&&&
 *[o][F-]{\bullet} \ar@{-}@(l,d)_{c} \ar@{-}[rr]_{b} && 
 *[o][F-]{\bullet} \ar@{-}@(r,d)^{a} 
 \\
 \ar@{}[u]|>>>>{\binom{1}{0}} && \ar@{}[u]|>>>>{\binom{0}{2}} && \ar@{}[u]|>>>>{\binom{2}{1}} &&&&
 &\ar@{}[ul]|>>>>{\binom{0}{0}} \ar@{}[ur]|>>>>{\binom{2}{2}} &
}
\]
\caption{$\textup{C}\Gamma_1$, the coupled Hanoi graph on 1 disk}
\label{f:CG1}
\end{figure}

Our goal is to provide some understanding of the coupled Hanoi graph $\textup{C}\Gamma_n$
on $n$ disks. The vertices of this graph are the $3^{2n}$ pairs of words
$\binom{u_T}{u_B}$ of length $n$ over $X$ (representing the top and the bottom
configuration on $n$ disks in a coupled configuration). Two vertices in
$\textup{C}\Gamma_n$ are connected by an edge labeled by $s$ in $\{a,b,c\}$ if the
coupled configurations represented by these vertices can be obtained from each
other by application of the move $s$. The coupled Hanoi graph on 1 disk is
depicted in Figure~\ref{f:CG1}.  
The coupled Hanoi graph $\textup{C}\Gamma_n$ is the Schreier graph, with respect to the
generating set $\{a,b,c\}$, of the action of $H$ on the pairs of words of length
$n$ in $X^*$ defined by   
\[
s \binom{u_T}{u_B} = \binom{s(u_T)}{s(u_B)}, 
\]
for $s$ in $\{a,b,c\}$.


\section{Twin Towers Switch}

In this section we provide an upper bound on the number of moves needed to solve
the Twin Towers Switch Problem. In the language of coupled Hanoi graphs the same
result is expressed as follows.  

\begin{theoremTTS'}
The distance between the coupled configurations
\[ \binom{000 \dots 0}{222 \dots 2} \qquad \text{ and } \qquad \binom{222 \dots 2}{000 \dots 0}\]
in the coupled Hanoi graph $\textup{C}\Gamma_n$ (on $n$ disks) is no greater than 
\[
a(n) =  
 \begin{cases}
   1, & n=1, \\
   \frac{4}{3} \cdot 2^n -\frac{(-1)^n}{3} , & n \geq 2.
 \end{cases} 
\]
\end{theoremTTS'}

\begin{proof}
Let $n=2$. We have, by using~\eqref{e:formula},  
\[
 ababa(00) = abab(10) = aba(12) = ab(02) = a(22) = 22.  
\] 
Since $ababa$ is a palindrome, it has order 2 (as a group element) and,
therefore, $ababa(22)=00$. Thus the distance between the initial and the final
coupled configurations is no greater than 5 (it can be shown that it is actually
5).  

Assume that $n$ is even and $n \geq 4$. Consider the sequence of
$a(n)=\frac{4}{3} \cdot 2^n - \frac{1}{3}$ moves  
\[ ababa(cacababa)^{\frac{1}{3}(2^{n-1}-2)}.\]
  
Notice the pattern 
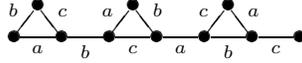
\[
\xymatrix@C=3pt@R=6pt{
 & 
 *[o][F-]{\bullet} \ar@{-}[dl]_{b} \ar@{-}[dr]^{c} &&&&
 *[o][F-]{\bullet} \ar@{-}[dl]_{a} \ar@{-}[dr]^{b} &&&& 
 *[o][F-]{\bullet} \ar@{-}[dl]_{c} \ar@{-}[dr]^{a} &&&&
 \\ 
 *[o][F-]{\bullet} \ar@{-}[rr]_{a} && *[o][F-]{\bullet} \ar@{-}[rr]_{b} &&
 *[o][F-]{\bullet} \ar@{-}[rr]_{c} && *[o][F-]{\bullet} \ar@{-}[rr]_{a} &&  
 *[o][F-]{\bullet} \ar@{-}[rr]_{b} && *[o][F-]{\bullet} \ar@{-}[rr]_{c} &&
 *[o][F-]{\bullet} 
}
\]
that repeats along the bottom edge in Figure~\ref{f:even-general}, indicating
that the result of the action of $cacababa$ and $(cba)^2$ on the leftmost vertex
in the pattern is the same and it is equal to the rightmost vertex in the
pattern (which is the leftmost vertex in the next occurrence of the pattern).
Therefore, by~\eqref{e:formula} and Lemma~\ref{l:classic}, 
\[ 
\begin{array}{rcl}
ababa(cacababa)^{\frac{1}{3}(2^{n-1}-2)}(000 \dots 0) 
  &= &ababa(cba)^{\frac{1}{3}(2^n-4)}(000 \dots 0) = 
  \\
  = ababaabc(cba)^{\frac{1}{3}(2^n-1)}(000 \dots 0) &= &
    abac(222 \dots 2) = 
  \\
  = aba(122\dots 2) &=& ab(022\dots 2) = \\
  = a(222 \dots 2) &=& 222\dots 2. 
\end{array}
\]

Since $ababa(cacababa)^{\frac{1}{3}(2^{n-1}-2)}$ is a palindrome it has order
$2$. Thus  
\[ ababa(cacababa)^{\frac{1}{3}(2^{n-1}-2)}(222 \dots 2) = 000 \dots 0 \]
and the distance between the initial and the final coupled configurations is no
greater than $a(n)$. 

Let $n=1$. The distance between the coupled configurations $\binom{0}{2}$ and
$\binom{2}{0}$ is 1 (see Figure~\ref{f:CG1}).  

Assume that $n$ is odd and $n \geq 3$. Consider the sequence of
$a(n)=\frac{4}{3} \cdot 2^n + \frac{1}{3}$ moves  
\[ aca(cbcbcaca)^{\frac{1}{3}(2^{n-1}-1)}.\]
  
Notice the pattern 
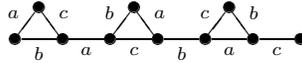
\[
\xymatrix@C=3pt@R=6pt{
 & 
 *[o][F-]{\bullet} \ar@{-}[dl]_{a} \ar@{-}[dr]^{c} &&&&
 *[o][F-]{\bullet} \ar@{-}[dl]_{b} \ar@{-}[dr]^{a} &&&& 
 *[o][F-]{\bullet} \ar@{-}[dl]_{c} \ar@{-}[dr]^{b} &&&&
 \\ 
 *[o][F-]{\bullet} \ar@{-}[rr]_{b} && *[o][F-]{\bullet} \ar@{-}[rr]_{a} &&
 *[o][F-]{\bullet} \ar@{-}[rr]_{c} && *[o][F-]{\bullet} \ar@{-}[rr]_{b} &&  
 *[o][F-]{\bullet} \ar@{-}[rr]_{a} && *[o][F-]{\bullet} \ar@{-}[rr]_{c} &&
 *[o][F-]{\bullet} 
}
\]
that repeats along the bottom edge in Figure~\ref{f:odd-general}, indicating
that the result of the action of $cbcbcaca$ and $(cab)^2$ on the leftmost vertex
in the pattern is the same and it is equal to the rightmost vertex in the
pattern (which is the leftmost vertex in the next occurrence of the pattern).
Therefore, by~\eqref{e:formula} and Lemma~\ref{l:classic}, 
\[ 
\begin{array}{rcl}
aca(cbcbcaca)^{\frac{1}{3}(2^{n-1}-1)}(000 \dots 0) 
  &= &aca(cab)^{\frac{1}{3}(2^n-2)}(000 \dots 0) = 
  \\
  = acabb(cab)^{\frac{1}{3}(2^n-2)}(000 \dots 0) &= &
    acab(222 \dots 2) = 
  \\
  = aca(022\dots 2) &=& ac(122\dots 2) = \\
  = a(222 \dots 2) &=& 222\dots 2.  
\end{array}
\]

Since $aca(cbcbcaca)^{\frac{1}{3}(2^{n-1}-1)}$ is a palindrome it has order $2$.
Thus  
\[ aca(cbcbcaca)^{\frac{1}{3}(2^{n-1}-1)}(222 \dots 2) = 000 \dots 0 \]
and the distance between the initial and the final coupled configurations is no
greater than $a(n)$. 
\end{proof}

\begin{remark}
There are several solutions of length $a(n)$, for $n \geq 2$. For instance,
another solution, for odd $n$, is  
\[ cac(ababacac)^{\frac{1}{3}(2^{n-1}-1)}\]
and, for even $n$, is 
\[ cbcbc(acacbcbc)^{\frac{1}{3}(2^{n-1}-2)}.\]
\end{remark}

We rephrase Conjecture~TTS as follows. 

\begin{conjectureTTS'}
The distance between the coupled configurations
\[ \binom{000 \dots 0}{222 \dots 2} \qquad \text{ and } \qquad \binom{222 \dots 2}{000 \dots 0}\]
in the coupled Hanoi graph $\textup{C}\Gamma_n$ (on $n$ disks) is equal to $a(n)$. 
\end{conjectureTTS'}


\section{Small Disk Shift}

For the considerations that follow, the concept of parity will be useful. 

\begin{definition}
For a configuration $u=x_1 \dots x_n$ in $X^*$ and $x \in X$, let $p_x(u)$ be
the parity of the number of appearances of the letter $x$ in $u$. For a coupled
configuration  $U = \binom{u_T}{u_B}$, let $p_x(U)$ be the parity of the sum of
the parities $p_x(u_T)$ and $p_x(u_B)$.  
\end{definition}

Call any of the configurations $0^n,1^n,2^n$ a \emph{corner configuration}. Call
a coupled configuration a \emph{corner coupled configuration} if at least one of
the configurations in it is a corner configuration. Application of $a_{ij}$ to
any non corner configuration changes the parities of both $i$ and $j$.
Therefore, application of $a_{ij}$ to a non corner coupled configuration does
not change any parities.  

\begin{theoremSDS'}[Small Disk Shift]
The distance between the coupled configurations 
\[ \binom{000 \dots 0}{100 \dots 0} \qquad \text{ and } \qquad \binom{100 \dots 0}{200 \dots 0}\]
in the coupled Hanoi graph $\textup{C}\Gamma_n$ (on $n$ disks) is 
\[
d(n) = 
 \begin{cases}
   2, & n=1, \\
   6, & n=2, \\
   2 \cdot 2^n, & n \geq 3.
 \end{cases} 
\]
\end{theoremSDS'}


\begin{proof}[Proof of \textup{Theorem~SDS$'$}: upper bound]
Assume that $n$ is even and $n \geq 4$. Consider the sequence of $2 \cdot 2^n$
moves  
\[ bab(abc)^{\frac{1}{3}(2^n-4)}a(cba)^{\frac{1}{3}(2^n-1)}c.\]
A simple and convincing way to verify that this sequence of moves accomplishes the goal would be to trace the action in Figure~\ref{f:even-general}. A more formal
approach, using~\eqref{e:formula} and Lemma~\ref{l:classic}, gives  
\[
\begin{array}{rcl}
bab(abc)^{\frac{1}{3}(2^n-4)}a(cba)^{\frac{1}{3}(2^n-1)}c(000 \dots 0) &=& 
    bab(abc)^{\frac{1}{3}(2^n-4)}a(cba)^{\frac{1}{3}(2^n-1)}(000 \dots 0) = \\
  = bab(abc)^{\frac{1}{3}(2^n-4)}a(222 \dots 2) &=& 
    bab(abc)^{\frac{1}{3}(2^n-4)}(222 \dots 2) = \\
  = babcba(abc)^{\frac{1}{3}(2^n-1)}(222 \dots 2) &=&
    babcba(000 \dots 0) = \\
  = babcb(100 \dots 0) &=& babc(120 \dots 0) = \\
  = bab(220 \dots 0) &=& ba(020 \dots 0) = \\
  = b(120 \dots 0) &=& 100 \dots 0,  
\end{array}
\]
and
\[
\begin{array}{rcl}
bab(abc)^{\frac{1}{3}(2^n-4)}a(cba)^{\frac{1}{3}(2^n-1)}c(100 \dots 0) &= 
   & bab(abc)^{\frac{1}{3}(2^n-4)}a(cba)^{\frac{1}{3}(2^n-1)}(200 \dots 0) = \\
  = bab(abc)^{\frac{1}{3}(2^n-4)}acba(200 \dots 0) &=& 
    bab(abc)^{\frac{1}{3}(2^n-4)}acb(210 \dots 0) = \\
  = bab(abc)^{\frac{1}{3}(2^n-4)}ac(010 \dots 0) &=& 
    bab(abc)^{\frac{1}{3}(2^n-4)}a(020 \dots 0) = \\
  = bab(abc)^{\frac{1}{3}(2^n-4)}(120 \dots 0) &=& bab(120 \dots 0) = \\
  = ba(100 \dots 0) &=& b(000 \dots 0) = \\
    &=& 200 \dots 0,  
\end{array}
\]
where, in the transition between the first and second row, we used the fact that
$\frac{1}{3}(2^n-1)$ is odd.  

Assume that $n$ is odd and $n \geq 3$. Consider the sequence of $2 \cdot 2^n$ moves 
\[ bab(abc)^{\frac{1}{3}(2^n-5)}bcba(cba)^{\frac{1}{3}(2^n-2)}.\]

A simple and convincing way to verify that this sequence of moves accomplishes the goal 
would be to trace the action in Figure~\ref{f:odd-general}. A more formal
approach, using~\eqref{e:formula} and Lemma~\ref{l:classic}, gives  
\[
\begin{array}{rcl}
bab(abc)^{\frac{1}{3}(2^n-5)}bcba(cba)^{\frac{1}{3}(2^n-2)}(000 \dots 0) &= & 
    bab(abc)^{\frac{1}{3}(2^n-5)}bcb(111 \dots 1) = \\
  = bab(abc)^{\frac{1}{3}(2^n-5)}aabc(111 \dots 1) &=& 
    baba(bca)^{\frac{1}{3}(2^n-5)}abc(111 \dots 1) = \\
  = baba(bca)^{\frac{1}{3}(2^n-5)}ab(211 \dots 1) &=& 
    baba(bca)^{\frac{1}{3}(2^n-5)}a(011 \dots 1) = \\
  = baba(bca)^{\frac{1}{3}(2^n-5)}(111 \dots 1) &=&
    babcbaabca(bca)^{\frac{1}{3}(2^n-5)}(111 \dots 1) = \\
  = babcbaa(bca)^{\frac{1}{3}(2^n-2)}(111 \dots 1) &=&
    babcba(000 \dots 0) = \\
  = babcb(100 \dots 0) &=& babc(120 \dots 0) = \\ 
  = bab(220 \dots 0) &=& ba(020 \dots 0) = \\ 
  = b(120 \dots 0) &=& 100 \dots 0,  
\end{array}
\]
and
\[
\begin{array}{rcl}
bab(abc)^{\frac{1}{3}(2^n-5)}bcba(cba)^{\frac{1}{3}(2^n-2)}(100 \dots 0) &= & 
    bab(abc)^{\frac{1}{3}(2^n-5)}b(100 \dots 0) = \\
  = bab(abc)^{\frac{1}{3}(2^n-5)}(120 \dots 0) &=& 
    bab(120 \dots 0) = \\
  = ba(100 \dots 0) &=& b(000 \dots 0) \\
    &=& 2000 \dots 0. 
\end{array} 
\]

When $n=1$, a solution of length 2 is given by the sequence of moves $ba$ and,
for $n=2$, a solution of length 6 is given by the  sequence of moves $bcacba$.  
\end{proof}

\begin{remark}
Note that the above sequences of moves of length $2 \cdot 2^n$ are not unique.
For instance, for even $n$, $n \geq 4$, one could use 
\[ caba(bac)^{\frac{1}{3}(2^n-4)}b(cab)^{\frac{1}{3}(2^n-1)}.\]
\end{remark}


\begin{proof}[Proof of \textup{Theorem~SDS$'$}: lower bound]
Since the 0-parities for the initial and final coupled configurations are 
\[ 
 p_0\binom{000 \dots 0}{100 \dots 0} = 1 \qquad \text{and} \qquad p_0\binom{100 \dots 0}{200 \dots 0} = 0, 
\]
somewhere on the way from the initial to the final coupled configuration the 
0-parity changes. This parity cannot be changed at the corner coupled
configurations $\binom{0^n}{v}$ and $\binom{v}{0^n}$, where $v$ is not a corner
configuration. Since the 0-parity must be changed, any sequence of moves that
starts at the initial coupled configuration $\binom{000\dots0}{100\dots0}$ and
accomplishes this change involves a corner $a$-loop of a corner $b$-loop
application in either the top or in the bottom configuration. The 4
possibilities are given as cases $\textup{Top}_a$, $\textup{Top}_b$, $\textup{Bot}_a$ and $\textup{Bot}_b$ (standing for top configuration involved in a corner $a$-loop, top configuration involved in a corner
$b$-loop, etc.) in Table~\ref{t:4forms}, 
\begin{table}[!ht]
\[ 
\xymatrix@R=12pt@C=15pt{
 \text{case} & \text{initial} & \ar[r]^{\text{0-parity}}_{\text{change}}& &  &
 \text{even }n& \text{odd }n 
 \\
 \textup{Top}_a: & \binom{00\dots0}{10\dots0} \ar[r] & \binom{2^n}{\ast} \ar[r]^{a} & 
 \binom{2^n}{\ast} \ar[r] & \binom{10\dots0}{\ast}, & {2 \cdot 2^n -2} 
 & {2 \cdot 2^n - 1} 
 \\ 
 \textup{Top}_b: & \binom{00\dots0}{10\dots0} \ar[r] & \binom{1^n}{\ast} \ar[r]^{b} & 
 \binom{1^n}{\ast} \ar[r] & \binom{10\dots0}{\ast}, & {2 \cdot 2^n -1}
 & {2 \cdot  2^n - 2} 
 \\    
 \textup{Bot}_a: & \binom{00\dots0}{10\dots0} \ar[r] & \binom{\ast}{2^n} \ar[r]_{a} & 
 \binom{\ast}{2^n} \ar[r] & \binom{\ast}{20\dots0}, & {2 \cdot 2^n -2} 
 & {2 \cdot 2^n - 2} 
 \\    
 \textup{Bot}_b: & \binom{00\dots0}{10\dots0} \ar[r] & \binom{\ast}{1^n} \ar[r]_{b} & 
 \binom{\ast}{1^n} \ar[r] & \binom{\ast}{20\dots0}, & { 2 \cdot 2^n -2} 
 & {2 \cdot 2^n - 2} 
} 
\]
\caption{Changing the 0-parity}
\label{t:4forms}
\end{table}
where, in each case, $\ast$ denotes some configuration different from the one
with which it is coupled. The last two columns provide the number of steps in
the unique shortest path of the given form, for even and odd number of disks.  

Note that the above considerations already show that $d(n)\geq 2 \cdot 2^n -2$
and that the largest disk has to be moved in at least one coupled set of disks.
 
Further, any element $g$ in $H$ for which
$g\binom{000\dots0}{100\dots0}=\binom{100\dots0}{200\dots0}$ must act on the
first letter as the permutation $(012)$, which is an even permutation.
Therefore, the length of $g$ must be even. To complete the proof, all we need to
show is that none of the shortest paths (sequences of moves) of length $2 \cdot 2^n-2$ implicitly
mentioned in Table~\ref{t:4forms} solves the Small Disk Shift Problem. 

For the unique shortest path $g$ of length $2 \cdot 2^n -2$ in Case $\textup{Top}_a$, even
$n$, such that for the top configuration we have $g(000\dots0)=100\dots0$,
tracing the action in Figure~\ref{f:even-general} for the bottom configuration,
we obtain  
\[ 
 (bca)^{\frac{1}{3}(2^n-1)}(cba)^{\frac{1}{3}(2^n-1)}(100\dots0) = 201\dots 1 \neq 200\dots0. \]

For the unique shortest path $g$ of length $2 \cdot 2^n -2$ in Case $\textup{Bot}_a$, even
$n$, such that for the bottom configuration we have $g(100\dots0)=200\dots0$,
tracing the action in Figure~\ref{f:even-general} for the top configuration, we
obtain  
\[ 
 (cbc)(abc)^{\frac{1}{3}(2^n-4)}(acb)^{\frac{1}{3}(2^n-1)}(000\dots0) = 101\dots 1 \neq 100\dots0.
\]

For the unique shortest path $g$ of length $2 \cdot 2^n -2$ in Case $\textup{Bot}_b$, even
$n$, such that for the bottom configuration we have $g(100\dots0)=200\dots0$,
tracing the action in Figure~\ref{f:even-general} for the top configuration, we
obtain  
\[ 
 ac(bac)^{\frac{1}{3}(2^n-4)}(bca)^{\frac{1}{3}(2^n-1)}c(000\dots0) = 102\dots 2 \neq 100\dots0.
\]

For the unique shortest path $g$ of length $2 \cdot 2^n -2$ in Case $\textup{Top}_b$, odd
$n$, such that for the top configuration we have $g(000\dots0)=100\dots0$,
tracing the action in Figure~\ref{f:odd-general} for the bottom configuration,
we obtain  
\[ 
 (bca)^{\frac{1}{3}(2^n-2)}ba(cba)^{\frac{1}{3}(2^n-2)}(100\dots0) = 222\dots 2 \neq 200\dots0.  
\]
 
For the unique shortest path $g$ of length $2 \cdot 2^n -2$ in Case $\textup{Bot}_a$, odd
$n$, such that for the bottom configuration we have $g(100\dots0)=200\dots0$,
tracing the action in Figure~\ref{f:odd-general} for the top configuration, we
obtain  
\[ 
 (acb)^{\frac{1}{3}(2^n-2)}a(bca)^{\frac{1}{3}(2^n-2)}c(000\dots0) = 111\dots 1 \neq 100\dots0. 
\]

Finally, for the unique shortest path $g$ of length $2 \cdot 2^n -2$ in Case
$\textup{Bot}_b$, odd $n$, such that for the bottom configuration we have
$g(100\dots0)=200\dots0$, tracing the action in Figure~\ref{f:odd-general} for
the top configuration, we obtain  
\[ 
 c(bca)^{\frac{1}{3}(2^n-2)}b(acb)^{\frac{1}{3}(2^n-2)}(000\dots0) = 122\dots 2 \neq 100\dots0. \qedhere
\]
\end{proof}


\section{General Problem}

In this section we describe the compatible coupled configurations (recovering
the result of D'Angeli and Donno from~\cite{d'angeli-d:hanoi}) and then provide
an upper bound on the distance between any compatible coupled configurations. 

In order to accomplish the goals of this section, we need a bit more information
on the Hanoi Towers group $H$. In particular, we rely on the self-similarity of
the action of $H$ on the tree $X^*$. More on self-similar actions in general can
be found in~\cite{nekrashevych:book-self-similar}. For our purposes the
following observations suffice.  

The action of $a$, $b$ and $c$ on $X^*$ given by~\eqref{e:formula} can be
rewritten in a recursive form as follows. For any word $u$ over $X$,  
\begin{alignat}{6}\label{e:action}
 &a(0u) &&=  1u, & \quad
 &b(0u) &&=  2u, & \quad
 &c(0u) &&=  0c(u), \notag
 \\
 &a(1u) &&=  0u, & 
 &b(1u) &&=  1b(u), &
 &c(1u) &&=  2u, 
 \\
 &a(2u) &&=  2a(u), &
 &b(2u) &&=  0u, &
 &c(2u) &&=  1u. \notag
\end{alignat}
This implies that, for any sequence $g$ of moves, there exist a permutation
$\pi_g$ of $X$ and three sequences of moves $g_0$, $g_1$ and $g_2$ such that,
for every word $u$ over $X$,   
\begin{equation}\label{e:action-general}
 g(0u) = \pi_g(0)g_0(u), \qquad g(1u) = \pi_g(1)g_1(u), \qquad g(2u) = \pi_g(2)g_2(u).
\end{equation}
The permutation $\pi(g)$ is called the \emph{root permutation} and it indicates the
action of $g$ on the first level of the tree (just below the root), while $g_0$,
$g_1$ and $g_2$ are called the \emph{sections} of $g$ and indicate the action of $g$
below the vertices on the first level. When~\eqref{e:action-general} holds, we
write  
\[
 g = \pi_g \ (g_0,g_1,g_2)  
\]
and call the expression on the right a \emph{decomposition} of $g$. Note
that~\eqref{e:action-general} may be correct for many different sequences of
moves $g_0$ (or $g_1$ or $g_2$), but all these sequences represent the same
element of the group $H$. Decompositions of the generators $a$, $b$ and $c$ are
given by  
\begin{equation}\label{e:decomposition}
 a = (01) \ (1,1,a), \qquad b = (02) \ (1,b,1),\qquad c = (12) \ (c,1,1),  
\end{equation}
where 1 denotes the empty sequence of moves (the trivial automorphism of the
tree). Two decompositions may be multiplied by using the formula
(see~\cite{nekrashevych:book-self-similar} or~\cite{grigorchuk-s:standrews}) 
\begin{equation}\label{e:product}
 gh = \pi_g \ (g_0,g_1,g_2) \ \pi_h \ (h_0,h_1,h_2)  = \pi_g\pi_h \ (g_{h(0)}h_0,g_{h(1)}h_1,g_{h(2)}h_2). 
\end{equation}
The decompositions of the generators $a$, $b$ and $c$ given
in~\eqref{e:decomposition} and the decomposition product
formula~\eqref{e:product} are sufficient to calculate a decomposition for any
sequence of moves. We refer to such calculations as decomposition calculations.

\begin{theorem}[D'Angeli and Donno~\cite{d'angeli-d:hanoi}]\label{t:d'angeli-d}
Two coupled configurations $U = \binom{u_T}{u_B}$ and $V = \binom{v_T}{v_B}$ on
$n$ disks are compatible if and only if the length of the longest common prefix
of $u_T$ and $u_B$ is the same as the length of the longest common prefix of
$v_T$ and $v_B$.  
\end{theorem}

\begin{remark}
Note that Theorem~\ref{t:d'angeli-d} implies that the $n+1$ sets 
$\textup{C}\Gamma_{n,0},\textup{C}\Gamma_{n,1},\dots,\textup{C}\Gamma_{n,n}$,  where $\textup{C}\Gamma_{n,i}$
consists of the coupled configurations $\binom{u_T}{u_B}$ such that the length of the longest common prefix of $u_T$ and $u_B$ is $i$, are the connected components of the coupled Hanoi graph $\textup{C}\Gamma_n$. The largest of these
sets is $\textup{C}\Gamma_{n,0}$. It consists of $6 \cdot 9^{n-1}$ vertices, which are the basic coupled configurations (defined in the introduction). More generally, the set $\textup{C}\Gamma_{n,i}$ has $3^i \cdot 6 \cdot 9^{n-1-i}$ vertices, for $i=0,\dots,n-1$, and $\textup{C}\Gamma_{n,n}$ has $3^n$ vertices (moreover, $\textup{C}\Gamma_{n,n}$ is canonically isomorphic to $\Gamma_n$ through the isomorphism $u \leftrightarrow \binom{u}{u}$). 

Since every tree automorphism preserves prefixes, the connected components of
the coupled Hanoi graph must be subsets of the sets $\textup{C}\Gamma_{n,i}$. Thus, only
the other direction (showing that each of the sets $\textup{C}\Gamma_{n,i}$ is connected)
is interesting and needs to be proved. 
\end{remark}

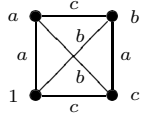
\begin{figure}[ht]
\[
\xymatrix{
 &
 *[o][F-]{\bullet} \ar@{}[l]|<<<{a} \ar@{-}[d]_{a} \ar@{-}[r]^{c} \ar@{-}[dr]^<<<<<<{b} & 
 *[o][F-]{\bullet} \ar@{}[r]|<<<{b} \ar@{-}[d]^{a}  &
 \\
 &
 *[o][F-]{\bullet} \ar@{}[l]|<<<{1} \ar@{-}[r]_{c} \ar@{-}[ur]_<<<<<<{b} & 
 *[o][F-]{\bullet} \ar@{}[r]|<<<{c} &
} 
\]
\caption{The Schreier graph of $A$ in $H$}
\label{f:schreier}
\end{figure}

Consider the subgroup $A=\langle cba,acb,bac\rangle \leq H$ (introduced
in~\cite{grigorchuk-n-s:oberwolfach1} and called Apollonian group, because its
limit space is the Apollonian gasket). It is known that this subgroup has index
4 in $H$ and $H/A=C_2 \times C_2$ (where $C_2$ is cyclic of order 2). A
sequence of moves $g$ belongs to $A$ if and only if the parities of the number
of occurrences of the moves $a$, $b$ and $c$ in $g$ are all odd or all even. The
elements $1,a,b,c$ form a transversal for $A$ in $H$. The
Schreier graph of the subgroup $A$ in $H$ is given in Figure~\ref{f:schreier}.
The vertices are denoted by the coset representatives (for instance, the vertex $b$ is the coset $bA$). 

\begin{lemma}
The Apollonian subgroup acts transitively on every level of the tree $X^*$. 
\end{lemma}

\begin{proof}
The claim follows from the fact that $H$ acts transitively on every level of the
tree and that, for every generator $s$ in $\{a,b,c\}$, there is a loop labeled
by $s$ in the Hanoi graph $\Gamma_n$.  

Indeed, if $g(u)=v$, for some sequence of moves $g$, and $g$ is in, say, the
coset $aA$, then $g'g(u)=v$ and $g'g$ is in $A$, where $g'=h^{-1}ah$ and $h$ is
any sequence of moves from $v$ to the vertex $2^n$ (note that $g' \in aA$ and $g'(v) = h^{-1}ah(v)=h^{-1}a(2^n)=h^{-1}(2^n)=v$).  
\end{proof}

\begin{remark}
A small modification of the above argument (using the corner loops to modify the
parity of the number of occurrences of any generator) shows that the commutator
subgroup $H'$ also acts transitively on every level of the tree. The fact that
$H'$ acts transitively was proved in a different way by D'Angeli and Donno and used in their proof of Theorem~\ref{t:d'angeli-d}. We provide a different proof of Theorem~\ref{t:d'angeli-d}, based on the transitivity of the action of $A$, enabling us to provide good estimates in the
General Problem for basic coupled configurations. 
\end{remark}

\begin{lemma}\label{l:large-connected}
The set $\textup{C}\Gamma_{n,0}$ of basic coupled configurations on $n$ disks is connected. 
\end{lemma}

\begin{proof}
Let $\binom{u_T}{u_B}$ and $\binom{v_T}{v_B}$ be coupled configurations in $\textup{C}\Gamma_{n,0}$. 

Since $H$ acts transitively on every level of the tree, there exists a sequence
of moves $h$ such that $h(u_T)=v_T$. Let $h(u_B)=v_B'$.  

Without loss of generality, assume that the top configuration $v_T$ starts by 2,
while the bottom configuration $v_B$ starts by 0. The configuration $v_B'$ may
start by either 0 or 1. If it starts by 1, a single application of the sequence
of 3 moves  
\[
 cab = (01) \ (a,cb,1), 
\]
does not affect $v_T$ (note the trivial section at 2), and changes the first letter in the bottom configuration
to 0. Thus, we may assume that both $v_B'$ and $v_B$ start by 0.  

We are interested in sequences of moves $g$ that do not affect any
configurations that start by 2 (and thus do not affect $v_T$) and keep the first
letter in the bottom configuration equal to 0. In other words, we are interested
in sequences of moves that decompose as  
\[
 g = (g_0,*,1), 
\]
where $*$ represents the section at 1, in which we are not interested. 

Three such sequences are (this can be verified by direct decomposition calculations)
\begin{align*}
 cabcab &= (cba,*,1) \\
 bacacaba &= (acb,*,1), \\ 
 bcbcacac &= (bac,*,1).  
\end{align*}
Since $\langle cba,acb,bac\rangle = A$, these three decompositions imply that,
for every sequence of moves $g_0$ in $A$, there exists a sequence of moves $g$
in $H$, and in fact in $A$, such that 
\[
 g = (g_0,*,1). 
\]

Let $v_B'=0v'$ and $v_B=0v$. Since $A$ acts transitively on each level of the
tree, there exists $g_0$ in $A$ such that $g_0(v')=v$. Therefore, there exists
$g$ in $A$ such that $g(v_B')=v_B$ and $g(v_T)=v_T$, completing the proof that
$\textup{C}\Gamma_{n,0}$ is connected.  
\end{proof}

The rest of the proof of Theorem~\ref{t:d'angeli-d} follows, essentially, the same steps as the original proof of D'Angeli and Donno and, being short, is included for completeness. Indeed, once it is known that the largest sets $\textup{C}\Gamma_{n,0}$ are connected, it is sufficient to observe that $H$ is a self-replicating group. 

\begin{lemma}
Hanoi Towers group $H$ is a self-replicating group of tree automorphisms, i.e., for every word $u$ over $X$ and every sequence of moves $g$ in $H$, there exists a sequence of moves $h$ in $H$ such that, for every word $w$ over $X$,  
\[
 h(uw) = ug(w).
\] 
\end{lemma}

\begin{proof}
Let $w$ be any word over $X$. Since 
\[
 a(2w) = 2a(w), \qquad cbc=2b(w), \qquad bcb(2w) = 2c(w), 
\]
it is clear that, for every sequence of moves $g$, there exists a sequence of moves $h$ such that $h(2w) = 2g(w)$. By symmetry, for every letter $x$ in $X$ and every sequence of moves $g$, there exists a sequence of moves $h$ such that 
\[
  h(xw) = xg(w) 
\]
and the claim easily extends to words over $X$ (and not just letters). 
\end{proof}

\begin{proof}[Proof of Theorem~\ref{t:d'angeli-d}]
Let $u$ and $u'$ be words of length $i$ and $\binom{uw_T}{uw_B}$ and $\binom{u'w_T'}{u'w_B'}$ be two coupled configurations in $\textup{C}\Gamma_{n,i}$. Since $H$ acts transitively on the levels of the tree, there exists a sequence of moves $h'$ in $H$ such that $h'\binom{uw_T}{uw_B} = \binom{u'w_T''}{u'w_B''}$, for some $w_T''$ and $w_B''$ (in fact, one may easily find such $h'$ for which $w_T''=w_T$ and $w_B''=w_B$, but this does not matter). Since $\textup{C}{\Gamma}_{n-i,0}$ is connected, there exists a sequence of moves $g$ such that $g\binom{w_T''}{w_B''} = \binom{w_T'}{w_B'}$. By the self-replicating property of $H$, there exists a sequence of moves $h$ in $H$ such that 
\[
 hh'\binom{uw_T}{uw_B} = h\binom{u'w_T''}{u'w_B''} = \binom{u'g(w_T'')}{u'g(w_B'')} = \binom{u'w_T'}{u'w_B'}. \qedhere
\]
\end{proof} 

\begin{theoremGP'}[General Problem for basic configurations]
The diameter $D(n)$ of the largest component $\textup{C}\Gamma_{n,0}$ of the coupled
Hanoi graph $\textup{C}\Gamma_n$ (on $n$ disks) satisfies, for $n \geq 3$, the
inequalities 
\[ 2 \times 2^n \leq D(n) \leq 3.\overline{66} \times 2^n. \]  
\end{theoremGP'}

\begin{proof}
We follow the proof of Lemma~\ref{l:large-connected}, but keep track of the lengths of the sequences of moves involved and, when we have a choice (and know how to make it), try to use short sequences. 

Let $U=\binom{u_T}{u_B}$ and $V=\binom{v_T}{v_B}$ be coupled configurations in
the largest component $\textup{C}\Gamma_{n,0}$ of the coupled Hanoi graph. Without loss
of generality, assume that the top configuration $v_T$ starts by 2, while the
bottom configuration $v_B$ starts by 0.  

There exists a sequence of moves $h$ of length at most $2^n+2$ such that
$h(u_T)=v_T$ and $h(u_B)= v_B'$, for some configuration $v_B'$ that starts by
0. Indeed, at most $2^n-1$ steps are needed to change the top configuration from
$u_T$ to $v_T$, and then at most three more steps (recall that $cab=(01)(a,cb,1)$) are
needed to make sure that the bottom configuration starts by 0.  

Let $v_B'=0v'$ and $v_B=0v$. We claim that there exists a sequence of moves
$g_0$ in $A$ such that $g_0(v')=v$ and the number of moves in the sequence $g_0$
is no greater than $2^n-1$. Indeed, if the shortest sequence of moves $g_s$
between $v'$ and $v$ happens to be in $A$ we may set $g_0=g_s$ (note that $v'$
and $v$ are vertices in the Hanoi graph $\Gamma_{n-1}$ of diameter $2^{n-1}-1$).
If $g_s$ happens to be, say, in the coset $aA$, we may set
$g_0=g^{(2)}ag^{(1)}$, where $g^{(1)}$ is the shortest sequence of moves from
$v'$ to $2^{n-1}$ and $g^{(2)}$ is the shortest sequence of moves from $2^{n-1}$
to $v$. The length of the sequence $g_0=g^{(2)}ag^{(1)}$ is no greater than
$2(2^{n-1}-1)+1= 2^n-1$. Since the sequence of moves $g_s^{-1}g^{(2)}g^{(1)}$
represents a closed path in the graph $\Gamma_{n-1}$ that does not go through
any of the corner loops and since all cycles in $\Gamma_{n-1}$ other than the
three corner loops are labeled by elements in $A$, the sequence
$g_s^{-1}g^{(2)}g^{(1)}$ is in $A$. Therefore  
\[
 g_0A = g^{(2)}ag^{(1)}A = ag^{(2)}g^{(1)}A = ag_sA = aaA = A, 
\]
which is what we needed. 

Direct decomposition calculations give
\begin{align*}
 bab(cba)^2bab &= (acaba,*,1) \\
 abc(acb)^2cba &= (babcb,*,1), \\ 
 cb(cba)^2bc &= (cbcac,*,1), 
\end{align*}
and therefore, for any $k \geq 0$, 
\begin{align}
 bab(cba)^{2k+2}bab &= (a(cab)^{k+1}a,*,1), \notag \\
 abc(acb)^{2k+2}cba &= (b(abc)^{k+1}b,*,1), \label{e:2k+2} \\ 
 cb(cba)^{2k+2}bc &= (c(bca)^{k+1}c,*,1). \notag 
\end{align}

This calculation justifies the entries in the top three rows of Table~\ref{t:move0}. 
\begin{table}[!ht]
\[
\begin{array}{r||r|r|r|r|r}
\text{case} & f & f_0 & \ell(f) & \ell(f_0) & \txt{ratio} \vspace{2mm} \\
\hline \hline
a \stackrel{-}\longleftarrow a & bab(cba)^{2k+2}bab & a(cab)^{k+1}a & 6k+12 & 3k+5 & 2.4 \\
\hline
b \stackrel{-}\longleftarrow b & abc(acb)^{2k+2}cba & b(abc)^{k+1}b & 6k+12 & 3k+5 & 2.4 \\
\hline
c \stackrel{-}\longleftarrow c & cb(cba)^{2k+2}bc & c(bca)^{k+1}c & 6k+10 & 3k+5 & 2 \\
\hline \hline 
 & cabcab & cba & 6 & 3 & 2 \\
c \stackrel{-}\longleftarrow a 
 & cabcb(cba)^{2k+2}bab & cb(cab)^{k+1}a & 6k+14 & 3k+6 & 2.34 \\
\hline 
 & bacacaba & acb & 8 & 3 & 2.67 \\
a \stackrel{-}\longleftarrow b  
 & bacacac(acb)^{2k+2}cba & ac(abc)^{k+1}b & 6k+ 16 & 3k+6 & 2.67 \\
\hline
 & bcbcacac & bac & 8 & 3 & 2.67 \\
b \stackrel{-}\longleftarrow c
 & bcbcacbcb(cba)^{2k+1}bc & ba(bca)^{k+1}c & 6k+14 & 3k+6 & 2.34 \\
\hline \hline
 & bcbcacbcab & baba & 10 & 4 & 2.5 \\
b \stackrel{-}\longleftarrow a
 & bcbcacbcb(cba)^{2k+2}bab & bab(cab)^{k+1}a & 6k+18 & 3k+7 & 2.58 \\
\hline
 & cabacaba & cbcb & 8 & 4 & 2 \\
c \stackrel{-}\longleftarrow b
 & cabacac(acb)^{2k+2}cba & cbc(abc)^{k+1}b & 6k + 16 & 3k+7 & 2.29 \\
\hline
 & babcbabc & acac & 8 & 4 & 2 \\
a \stackrel{-}\longleftarrow c
 & bab(cba)^{2k+3}bc & aca(bca)^{k+1}c & 6k+14 & 3k+7 & 2 \\
\hline \hline 
\end{array}
\]
\caption{Sequences of moves fixing $v_T$ and moving $v_B'$}
\label{t:move0}
\end{table}
In this table, $f$ is a sequence of moves and $f_0$ is the corresponding  section at 0. The first letter of any word is fixed by $f$ and the section at 2 is trivial. In other words, $f$ decomposes as 
\[
 f = (f_0,*,1). 
\]
The lengths of the sequences $f$ and $f_0$, as written, are $\ell(f)$ and $\ell(f_0)$, and the ratio in the last column is the ratio $\ell(f)/\ell(f_0)$ (in the rows that depend on $k$, the ratio is the maximum possible ratio, taken for $k \geq 0$ and rounded up).  

The entries in the remaining rows in Table~\ref{t:move0} are easy to verify. For
instance, for case $c \stackrel{-}\longleftarrow a$, by direct decomposition calculation,  
\begin{equation}\label{e:cba}
 cabcab = (cba,*,1)
\end{equation}
and the entry in the next row is obtained simply by multiplying the
equality~\eqref{e:cba} and the first equality in~\eqref{e:2k+2} 
\begin{align*}
 cabcb(cba)^{2k+2}bab &= (cabc)(abba)b(cba)^{2k+2}bab) = (cabcab)(bab(cba)^{2k+2}bab) = \\
  &= (cba,*,1)(a(cab)^{k+1}a,*,1) = (cbaa(cab)^{k+1}a,*,1) = \\
  &= (cb(cab)^{k+1}a,*,1).
\end{align*}
All other cases are equally easy to verify (by verifying directly the basic
case, and then multiplying it by a corresponding equality from~\eqref{e:2k+2} to
obtain the cases depending on $k$). 

Consider $g_0$ as defined above. There is no occurrence of $aa$, $bb$ or $cc$ in
this sequence (since we always chose the shortest paths as we built $g_0$) and
it is in $A$. The sequence $g_0$ is a product of factors each of which has the
form of one of the entries in column $f_0$ in Table~\ref{t:move0} or their
inverses. Moreover, the decomposition is such that the length of $g_0$ is the
sum of the lengths of the factors. Indeed, the entries in column $f_0$ and their
inverses are all possible sequences of moves in $A$ without occurrence of $aa$,
$bb$ or $cc$ for which no proper suffix is in $A$. Such sequences correspond
precisely to paths without backtracking in the Schreier graph in
Figure~\ref{f:schreier} that start at 1, end at 1 and do not visit the vertex 1
except at the very beginning and at the very end. There are 18 such types of
paths, three choices for the first step ($a$, $b$ or $c$) to leave vertex 1,
three choices for the last step ($a$, $b$ or $c$) to go back to vertex 1, and
two choices for the orientation (order) used to loop around the three vertices (cosets) 
$a$, $b$ and $c$ before the return to 1 (positive or negative orientation). The column $f_0$ in the table only
lists the 9 possible paths with negative orientation (and classifies the 9
cases by the first and last move), since the other 9 are just
inverses of the entries in the table. For instance, the notation $c \stackrel{-}{\longleftarrow} a$ indicates paths (sequences of moves) that start by the move $a$ and end by the move $c$. 

Once $g_0$ is appropriately factored, Table~\ref{t:move0} can be used to define
$g$ of length no greater than $2.\overline{66} \ell(g_0) \leq 2.\overline{66}
(2^n-1)$ such that $g\binom{v_T}{v_B'}=\binom{v_T}{v_B}$.  

Thus, we may arrive from the initial coupled configuration $\binom{u_T}{u_B}$ to
the final coupled configuration $\binom{v_T}{v_B}$ in no more than $(2^n+2) + 2.\overline{66} (2^n-1) \leq 3.\overline{66} \times 2^n$ moves. 
\end{proof}

It is evident that good understanding of the structure of $\textup{C}\Gamma_{n,0}$, for all $n$, provides good understanding of $\textup{C}\Gamma_{n,i}$, for all $n$ and $i$. For instance, the understanding of the graphs $\textup{C}\Gamma_{1,0}$ (6 vertices, diameter 2) and $\textup{C}\Gamma_{2,0}$ (54 vertices, diameter 6) enabled the author to determine the exact values of the diameter of the two smallest nontrivial components $\textup{C}\Gamma_{n,n-1}$ and $\textup{C}\Gamma_{n,n-2}$, for any number of disks. For instance, the diameter of
$\textup{C}\Gamma_{n,n-1}$ is, for $n \geq 1$, equal to 
\[ 
 \frac{7}{6} 2^n - \frac{3+(-1)^n}{6}. 
\]
The details will appear in a future work. 

\def\cprime{$'$}


\end{document}